\documentclass[11pt]{article}
\usepackage[T1]{fontenc}
\usepackage[utf8]{inputenc}
\usepackage[margin=1in]{geometry}
\usepackage{graphicx}
\usepackage{amsmath,amssymb}
\usepackage{authblk}
\usepackage{cite}

\title{Optimal Control of Pandemic Dynamics via Model Predictive Control:\\
A Health--Economic Trade-off Analysis}

\author[1]{Lokman Rachid Melhani\thanks{Corresponding author:
\texttt{lokmanrachid.melhani@unipa.it}}}
\author[2]{Lars Gr\"une}
\author[1]{Antonino Sferlazza}
\author[1]{Dominique Persano Adorno}
\author[1]{Filippo D'Ippolito}
\author[3]{Alberto Firenze}

\affil[1]{Department of Engineering, University of Palermo, Palermo, Italy}
\affil[2]{Department of Mathematics, University of Bayreuth, Bayreuth, Germany}
\affil[3]{Department of Internal Medicine ``Promise'', University of Palermo,
90127 Palermo, Italy}

\date{}

\begin{document}
\maketitle

\begin{abstract}
This paper addresses the optimal control of epidemic dynamics under conflicting
socio-economic objectives. We propose an economic Model Predictive Control (MPC)
framework, applied to an extended SEIR--V
(Susceptible--Exposed--Infected--Recovered--Vaccinated) compartmental model to
govern the spread of an infectious disease while minimizing economic disruption.
The control problem is formulated as a constrained nonlinear optimization
problem, in which the controller dynamically adjusts social interaction levels
(transmission rate $\beta$) and vaccination efforts to minimize a composite cost
function that penalizes fatalities, healthcare capacity violations, and economic
losses. We conduct a rigorous sensitivity analysis of the prediction horizon
$N$, demonstrating that the closed loop is robust to the horizon choice and that
$N=35$ days minimizes the realized cost. Furthermore, both the closed-loop
solution and an open-loop turnpike analysis across diverse initial conditions
reveal that the celebrated ``Hammer and Dance'' mitigation strategy emerges
naturally as the mathematical optimum: the optimal trajectories anchor to a
unique suppression turnpike (maximum lockdown) to drive hospitalizations toward
the disease-free equilibrium before progressively reopening the economy. Through
a turnpike-based argument we establish practical asymptotic stability of the
optimal operating point, providing a mathematically grounded decision-support
tool for pandemic policy.
\end{abstract}


\bigskip

\section{Introduction}\label{sec:intro}
The outbreak of the COVID-19 pandemic exposed a fundamental vulnerability in
global crisis management: the lack of quantitative, feedback-driven tools to
balance competing societal objectives. Governments worldwide were forced to
navigate a ``lives versus livelihoods'' dilemma, choosing between strict
Non-Pharmaceutical Interventions (NPIs) such as lockdowns, which suppress viral
transmission but cripple economic activity, and relaxed measures that preserve
economic flow at the cost of public health surges.

From a systems engineering perspective, an epidemic is a dynamic, nonlinear
process governed by biological parameters such as incubation periods and
transmission rates. Traditional epidemiological approaches
\cite{kermack1927contribution} have largely focused on \textit{prediction},
forecasting infection curves under static assumptions. However, effective
pandemic management is inherently a \textit{control problem}. The objective is
not merely to observe the system but to regulate it: steering the state
trajectory (infection levels) toward a safe equilibrium while minimizing the
control effort (economic cost).

This paper addresses the optimal control of pandemic dynamics using
\textbf{Model Predictive Control (MPC)}. Unlike open-loop strategies, which fix
policies months in advance and fail to adapt to stochastic disturbances, MPC
solves a finite-horizon optimization problem at every time
step.\cite{camacho2007constrained} This approach allows for real-time adaptation
to changing infection rates and ensures that distinct system constraints are
rigorously respected, as demonstrated in recent applications of control theory
to COVID-19.\cite{morato2020optimal,kohler2020robust}

\subsection{Contributions}\label{subsec:contrib}
The contributions of this work are threefold:
\begin{enumerate}
\item \textbf{Integrated and calibrated Modeling.} We extend the classical SEIR
framework to an \textbf{SEIR--V} model, incorporating vaccination as a secondary
control input alongside social distancing. This model was calibrated to capture
the dynamics of the COVID-19 pandemic in Italy during the critical vaccination
phase (January--May 2021); the explicit mathematical calibration and parameter
identification using Italian healthcare data are detailed in a companion
manuscript.\cite{melhani2026nine}
\item \textbf{Comprehensive Multi-Objective Cost Formulation.} We formulate a
composite objective function that explicitly quantifies the trade-off between
health system saturation, economic loss, and vaccination costs, using a deadzone
soft constraint to encode finite intensive-care capacity.
\item \textbf{Turnpike-Based Stability Analysis.} Within the framework of
economic MPC \emph{without} terminal constraints, we establish -- through the
turnpike property -- the practical asymptotic stability of the optimal operating
point, characterize that point as an open-loop unstable disease-free
equilibrium, and show that the resulting optimal policy is the Hammer-and-Dance
strategy.
\end{enumerate}
This framework provides a mathematically grounded decision-support tool that may
help pandemic policy move from reactive heuristics to optimal, proactive
regulation.

\section{Problem Formulation and Pandemic Modeling}\label{sec:model}
We consider the dynamics of an infectious disease propagating through a
population of size $N$. The model structure is based on the compartmental
framework described in \cite{anderson1991infectious}, extended here to capture
the heterogeneity of infection severity, the impact of vaccination, and waning
immunity. We employ an extended \textbf{SEIR--V} compartmental model.

\subsection{State Space Representation}\label{subsec:state}
The system state at time $t$ is defined by the vector $x(t)\in\mathbb{R}^{9}$:
\begin{equation}
x(t) = [\,S,\; E,\; V,\; I_1,\; I_2,\; P,\; H,\; R,\; F\,]^{\top}
\end{equation}
The compartments are defined as follows:
\begin{itemize}
\item $S(t)$: Susceptible individuals
\item $E(t)$: Exposed individuals (infected but not yet infectious)
\item $V(t)$: Vaccinated individuals
\item $I_1(t)$: Infected with the original virus
\item $I_2(t)$: Infected with the mutant virus
\item $P(t)$: Super-Spreaders
\item $H(t)$: Hospitalized individuals
\item $R(t)$: Recovered individuals
\item $F(t)$: Fatalities (Deceased)
\end{itemize}

\subsection{System Dynamics}\label{subsec:dyn}
The time-evolution of the epidemic is governed by the following system of
nonlinear ordinary differential equations. The model accounts for distinct
transmission rates for the infectious classes ($\beta,\beta',\beta_2$), waning
natural immunity ($\delta$), and waning vaccine efficacy ($\psi$). For
compactness we write the aggregate force of infection as
\begin{equation}
\lambda(t)=\beta(t)\frac{I_1}{N}+\beta'(t)\frac{P}{N}+\beta_2(t)\frac{I_2}{N}.
\label{eq:foi}
\end{equation}
\begin{footnotesize}
\begin{align}
\frac{dS}{dt} &= \Lambda + \psi V + \delta R - \lambda(t) S
 - (\mu + w_{1}(t)) S \label{eq:dS}\\
\frac{dE}{dt} &= \lambda(t) S + (1-\sigma)\,\lambda(t)\, V
 - (\mu + \kappa) E \label{eq:dE}\\
\frac{dV}{dt} &= w_{1}(t) S - (1-\sigma)\,\lambda(t)\, V
 - (\mu + \psi) V \label{eq:dV}\\
\frac{dI_{1}}{dt} &= \kappa \rho_{1} E
 - (\gamma_{a}+\gamma_{i}+\delta_{i}+m+r_{1}+\mu) I_{1} \label{eq:dI1}\\
\frac{dI_{2}}{dt} &= \kappa (1-\rho_{1}-\rho_{2}) E + m I_{1}
 - (\gamma_{a}+\gamma_{i}+\delta_{i}+r_{2}+\mu) I_{2} \label{eq:dI2}\\
\frac{dP}{dt} &= \kappa \rho_{2} E
 - (\gamma_{a}+\gamma_{i}+\delta_{p}+\mu) P \label{eq:dP}\\
\frac{dH}{dt} &= \gamma_{a}(I_{1}+I_{2}+P)
 - (\gamma_{r}+\delta_{h}+\mu) H \label{eq:dH}\\
\frac{dR}{dt} &= \gamma_{i}(I_{1}+I_{2}+P) + \gamma_{r} H
 + r_{1} I_{1} + r_{2} I_{2} - (\mu+\delta) R \label{eq:dR}\\
\frac{dF}{dt} &= \delta_{i}(I_{1}+I_{2}) + \delta_{p} P
 + \delta_{h} H \label{eq:dF}
\end{align}
\end{footnotesize}
The breakthrough-infection term in \eqref{eq:dE}--\eqref{eq:dV} carries the
\emph{full} force of infection \eqref{eq:foi}, including the super-spreader
pathway $\beta'(t)P/N$, consistent with the calibrated
model.\cite{melhani2026nine}

\subsection{Model Illustration}\label{subsec:illus}
This section details the biological and epidemiological significance of each
equation.

\textbf{Susceptible (S):} Population inflow occurs through new births ($\Lambda$)
and immunity waning from both vaccinated ($\psi V$) and recovered ($\delta R$)
individuals. Outflow results from three infection pathways (original strain,
mutant strain, super-spreaders) through the time-varying rates
$\beta(t),\beta'(t),\beta_2(t)$, along with natural mortality ($\mu$) and
vaccination uptake ($w_1(t)S$).

\textbf{Exposed (E):} Inflow occurs through infections of susceptibles from the
three transmission sources and breakthrough infections in vaccinated
individuals, where the $(1-\sigma)$ factor encodes reduced susceptibility due to
vaccine protection. Outflow is governed by progression to infectious states
($\kappa E$).

\textbf{Infected, original strain ($I_1$):} Outflow pathways are hospitalization
($\gamma_a I_1$), community recovery ($\gamma_i I_1$), strain-specific direct
recovery ($r_1 I_1$), disease-induced mortality ($\delta_i I_1$), and mutation to
the new strain ($m I_1$).

\textbf{Infected, mutant strain ($I_2$):} Inflow arrives through both direct
infection from exposed individuals and mutation from the original strain
($m I_1$).

\textbf{Vaccinated (V):} Inflow from susceptibles via vaccination ($w_1(t)S$);
outflow through natural mortality, breakthrough infection, and waning vaccine
immunity ($\psi V$).

\textbf{Super-spreaders (P):} Inflow is the fraction ($\rho_2$) of exposed
individuals developing enhanced transmission potential.

\textbf{Hospitalized (H):} Outflow through recovery with treatment
($\gamma_r H$), in-hospital mortality ($\delta_h H$), and natural mortality.

\textbf{Recovered (R):} Inflow through the several recovery pathways; outflow
through natural mortality and waning natural immunity ($\delta R$).

\textbf{Fatalities (F):} Accumulates disease-induced deaths from the infectious
communities ($\delta_i(I_1+I_2)$), super-spreaders ($\delta_p P$), and
hospitalized patients ($\delta_h H$).

\subsection{Control Inputs}\label{subsec:inputs}
The system is actuated by two control variables $u(t)=[\beta(t),w_1(t)]^{\top}$:
\begin{enumerate}
\item \textbf{Transmission Rate ($\beta(t)$):} the level of social interaction
and a proxy for economic openness. The associated rates $\beta'(t)$ and
$\beta_2(t)$ scale with $\beta(t)$ to reflect the relative transmissibility of
super-spreaders ($P$) and of the mutant strain ($I_2$):
\begin{equation}
\beta'(t)=c_P\,\beta(t),\qquad \beta_2(t)=c_2\,\beta(t),\qquad c_P=c_2=1.5,
\label{eq:multipliers}
\end{equation}
with the multipliers determined by the calibration.\cite{melhani2026nine}
\item \textbf{Vaccination Rate ($w_1(t)$):} determines the rate at which
susceptibles are moved to the vaccinated compartment $V$, subject to logistical
constraints.
\end{enumerate}

\section{The Optimal Control Problem}\label{sec:ocp}
The management of a pandemic is fundamentally a multi-objective optimal control
problem subject to severe state constraints and competing socio-economic
priorities. Compared to the standard SIR or basic SEIR models frequently used in
early COVID-19 control literature, our extended SEIR--V framework provides a
higher-fidelity representation of the epidemic's complex later stages. By
explicitly incorporating vaccination dynamics, super-spreaders, and mutant
strains, the model allows the controller to coordinate pharmaceutical
interventions alongside non-pharmaceutical ones.

The goal is to determine a control law $u^{*}(t)=[\beta(t),w_1(t)]^{\top}$ that
minimizes the cumulative societal cost over a finite horizon. Unlike standard
setpoint-tracking problems, this multi-objective pandemic control balances
inherently conflicting objectives: minimizing infection peaks and fatalities
while simultaneously maximizing economic activity and managing the logistical
effort of vaccination.

\subsection{Discrete-Time Formulation}\label{subsec:disc}
While MPC can be formulated in continuous time, its practical implementation via
numerical optimization requires a discrete-time framework. We therefore
discretize the continuous dynamics \eqref{eq:dS}--\eqref{eq:dF} with a sampling
time $T_s=1$~day. Let $x_k$ denote the state at day $k$ and $u_k$ the control
held constant on $[k,k+1)$. The discretized dynamics are
\begin{equation}
x_{k+1}=f(x_k,u_k),\label{eq:disc}
\end{equation}
where $f(\cdot)$ is obtained by fourth-order Runge--Kutta integration of the
continuous ODEs. This transforms the infinite-dimensional continuous problem
into a finite-dimensional nonlinear program (NLP), consistent with standard MPC
theory.\cite{rawlings2017model,grune2016nonlinear}

\subsection{The Multi-Objective Cost Function}\label{subsec:cost}
We formulate the cost $J_N$ over a prediction horizon of $N$ days as a weighted
sum of five quadratic terms,
\begin{equation}
J_N(x_0,\mathbf{u})=\sum_{k=0}^{N-1}\ell(x_k,u_k),\quad
\ell = J_{\mathrm{econ}}+J_{\mathrm{health}}+J_{\mathrm{death}}
+J_{\mathrm{vax}}+J_{\mathrm{smooth}}.
\end{equation}

\subsubsection{Economic Cost ($J_{\mathrm{econ}}$)}
The relationship between disease transmission and economic activity is well
established.\cite{acemoglu2021optimal} We model economic loss as a quadratic
deviation from the ``Business-as-Usual'' transmission rate $\beta_{\max}$:
\begin{equation}
J_{\mathrm{econ}}(k)=w_{\mathrm{econ}}
\left(\frac{\beta_{\max}-\beta_k}{\beta_{\max}-\beta_{\min}}\right)^{2}.
\end{equation}
The marginal cost of restrictions thus increases nonlinearly: mild distancing is
cheap, whereas full lockdown is disproportionately costly.

\subsubsection{Healthcare Capacity Violation ($J_{\mathrm{health}}$)}
To prevent collapse of the healthcare system, with finite ICU capacity
$H_{\mathrm{safe}}$, we use a \textbf{soft constraint} barrier
penalty:\cite{kerrigan2000soft}
\begin{equation}
J_{\mathrm{health}}(k)=w_{\mathrm{health}}
\left(\max\!\Big\{0,\frac{H_k-H_{\mathrm{safe}}}{H_{\mathrm{scale}}}\Big\}\right)^{2}.
\end{equation}
This term is zero while $H_k\le H_{\mathrm{safe}}$ and rises quadratically beyond
it, allowing temporary violation if physically unavoidable while exerting immense
pressure to return to safe levels.

\subsubsection{Fatality Minimization ($J_{\mathrm{death}}$)}
To value human life beyond economics, we penalize the accumulation of fatalities
through their daily increment,
\begin{equation}
J_{\mathrm{death}}(k)=w_{\mathrm{death}}\,(F_{k+1}-F_k)^{2},
\end{equation}
thereby targeting peak mortality periods. Since $F$ is a cumulative fatality
counter, its daily increment $F_{k+1}-F_k\ge 0$ is the (non-negative) number of
deaths on day $k$; the quadratic form therefore weighs days of concentrated
mortality more heavily than the same number of deaths spread over time. Because
this increment cannot be negative in the model, the squared penalty never acts on
a reduction of the cumulative count.

\subsubsection{Vaccination Effort ($J_{\mathrm{vax}}$)}
Mass vaccination faces logistic and financial constraints, modeled as
\begin{equation}
J_{\mathrm{vax}}(k)=w_{\mathrm{vax}}\left(\frac{w_{1,k}}{w_{\max}}\right)^{2},
\end{equation}
reflecting the increasing marginal cost of rapidly reaching the whole population.

\subsubsection{Policy Smoothing ($J_{\mathrm{smooth}}$)}
Social adherence depends on policy consistency. Without this penalty the optimum
exhibits ``bang-bang'' chattering of $\beta(t)$ between its bounds, which causes
non-compliance and is practically unenforceable. We therefore penalize the rate
of change of $\beta$:
\begin{equation}
J_{\mathrm{smooth}}(k)=w_{\mathrm{smooth}}\,(\beta_k-\beta_{k-1})^{2}.
\end{equation}

\subsection{The Constrained Optimization Problem}\label{subsec:nlp}
The optimal control problem (OCP) solved at each time step $t$ is
\begin{align}
\min_{\mathbf{u}}\quad & \sum_{k=0}^{N-1}\ell(x_k,u_k) \label{eq:ocp}\\
\text{s.t.}\quad & x_{k+1}=f(x_k,u_k),\quad k=0,\dots,N-1 \nonumber\\
& x_0=x(t)\quad(\text{current state feedback}) \nonumber\\
& \beta_{\min}\le\beta_k\le\beta_{\max} \nonumber\\
& 0\le w_{1,k}\le w_{\max} \nonumber
\end{align}
optimized over $\mathbf{u}=\{u_0,\dots,u_{N-1}\}$. The initial state $x_0$ is
updated at every step from the current epidemic state; in deployment this state
is reconstructed from partial measurements by the Extended Kalman Filter
developed for this model in a companion work,\cite{melhani2026extended} while the
closed-loop results reported here assume full-state feedback.

\section{Model Predictive Control Strategy}\label{sec:mpc}
The dynamic, uncertain, nonlinear nature of an epidemic renders classical
open-loop optimal control ineffective: a policy computed at $t=0$ inevitably
diverges from the true optimum as reality unfolds. We therefore implement a
\textbf{Receding Horizon Control (RHC)} strategy, introducing feedback by
re-optimizing at every sampling instant.

\subsection{The Receding Horizon Principle}\label{subsec:rhc}
At each instant $k$ the controller solves the OCP over the prediction horizon
$N$,\cite{rawlings2017model}
\begin{equation}
\mathbf{u}^{*}(k)=\arg\min_{\mathbf{u}} J_N(x(k),\mathbf{u}),
\end{equation}
yielding an optimal sequence $\mathbf{u}^{*}(k)=\{u^{*}(0\mid k),\dots,
u^{*}(N-1\mid k)\}$. Only the first element $u_{\mathrm{mpc}}(k)=u^{*}(0\mid k)$
is applied; at $k+1$ the horizon shifts and the optimization repeats with the new
measured state $x(k+1)$. This closed-loop feedback confers inherent
robustness:\cite[Chapter~7]{grune2016nonlinear} rather than following an outdated
prediction, the controller continuously adjusts the economic openness $\beta(t)$
to compensate for real-world deviations. We emphasize that the robustness
addressed in this work is that of receding-horizon feedback to initial-state
perturbations and to the residual mismatch absorbed by re-optimization; an
explicit treatment of parametric and measurement uncertainty (e.g. via stochastic
or tube-based MPC) is deferred to future work.

\subsection{Numerical Implementation}\label{subsec:impl}
The OCP is a \textbf{Nonlinear Programming (NLP)} problem. We employ a direct
single-shooting approach combined with a \textbf{Sequential Quadratic Programming
(SQP)} algorithm.\cite{biegler2010nonlinear}
\begin{enumerate}
\item \textbf{Solver.} We use the \texttt{fmincon} solver in MATLAB; the SQP
algorithm approximates the Lagrangian by a quadratic model and the nonlinear
constraints by linearizations, offering superlinear local
convergence.\cite{nocedal2006numerical}
\item \textbf{Integration.} The continuous dynamics are discretized by a
\textbf{fourth-order Runge--Kutta} method. To justify the fixed step $T_s=1$~day,
we refined the integration grid to $T_s=0.5$ and $0.25$~days and found the
controlled trajectories -- in particular the peak hospitalization -- essentially
unchanged, confirming that a one-day step ensures high-fidelity prediction while
keeping the NLP tractable.
\item \textbf{Warm Starting.} The solver is warm-started at each step by shifting
the previous optimal solution, an essential device for real-time nonlinear
MPC.\cite{diehl2002real}
\end{enumerate}

\subsection{Feasibility, Turnpike and Stability}\label{subsec:feas}
A common challenge in MPC is recursive feasibility. First, the healthcare
constraint $H\le H_{\mathrm{safe}}$ is a \textit{soft constraint} via the penalty
of Section~\ref{subsec:cost};\cite{kerrigan2000soft} this keeps the optimization
feasible even when the system starts in a violated state, the solver recovering
safety as cheaply as possible. Second, regarding stability, we employ neither a
terminal constraint set nor a terminal penalty -- both are typically intractable
for high-dimensional nonlinear epidemic models. Instead we rely on the framework
of Economic MPC without terminal constraints.\cite{grune2016nonlinear}

\paragraph{Characterization of the optimal operating point.}
We first identify the steady state about which the economic MPC operates. Let
$\mathcal{M}_0=\{x:\;E=I_1=I_2=P=H=0\}$ denote the disease-free manifold. On
$\mathcal{M}_0$ the force of infection \eqref{eq:foi} vanishes
($\lambda\equiv0$); consequently the infected compartments remain identically
zero -- $\mathcal{M}_0$ is forward invariant -- and the transmission rate $\beta$
no longer influences the dynamics, so that maximal openness $\beta=\beta_{\max}$
is admissible there at no epidemiological cost. Within $\mathcal{M}_0$ the
residual demographic dynamics admit, for a constant vaccination rate $w_1$, the
equilibrium
\begin{equation}
S^{\star}=\frac{\Lambda(\mu+\psi)}{\mu(\mu+\psi+w_1)},\qquad
V^{\star}=\frac{w_1\,\Lambda}{\mu(\mu+\psi+w_1)},\qquad R^{\star}=0,
\label{eq:dfe}
\end{equation}
with the infected and hospitalized states zero and the fatality counter $F$
constant. Since every stage-cost term is non-negative and they vanish
simultaneously at $\beta=\beta_{\max},\,w_1=0$, the optimal operating point
$(x^{\star},u^{\star})$ is the disease-free state with $u^{\star}=(\beta_{\max},0)$
and $\ell(x^{\star},u^{\star})=0$; it is a steady state,
$f(x^{\star},u^{\star})=x^{\star}$.

Two features make this point nontrivial and explain why open-loop optimization is
inadequate. First, on the epidemic timescale of months the demographic rates are
negligible ($\Lambda,\mu=\mathcal{O}(10^{-5})$ per day), so $\mathcal{M}_0$
itself -- rather than its far-field demographic equilibrium \eqref{eq:dfe} -- is
the operative turnpike target. Second, and decisively, $\mathcal{M}_0$ is
\emph{open-loop unstable} in the transversal (infected) directions whenever the
basic reproduction number at maximal openness exceeds unity,
$R_0(\beta_{\max})>1$, as holds for the calibrated Alpha-variant parameters: an
arbitrarily small infection seed grows exponentially under the constant control
$u^{\star}$. The optimal steady state is therefore an \emph{unstable} equilibrium
that the economic MPC feedback must actively stabilize -- precisely the regime in
which the turnpike property, rather than open-loop convergence, governs the
stability of the closed loop.

\noindent\textbf{Definition 1 (Turnpike property).} \textit{The OCP
\eqref{eq:ocp} has the turnpike property at $x^{\star}$ if optimal trajectories
spend all but a bounded, horizon-independent number of steps within any
$\varepsilon$-neighborhood of $x^{\star}$, departing in an initial approach
phase and a terminal exit phase.}

For economic MPC without terminal constraints, the turnpike property is the
structural mechanism through which closed-loop stability is obtained: when the
optimal control problem exhibits a turnpike at $x^{\star}$ and the system is
suitably controllable, a sufficiently long horizon $N$ renders $x^{\star}$
practically asymptotically stable, and state constraints are accommodated without
stabilizing terminal
sets.\cite{grune2013economic,grune2014asymptotic,faulwasser2018economic,boccia2014stability}
We summarize the expected behaviour informally, deferring the rigorous
statement to the literature.

\noindent\textbf{Remark 1.} When the OCP \eqref{eq:ocp} exhibits the turnpike
property at $(x^{\star},u^{\star})$ and the system is suitably controllable,
economic MPC without terminal constraints is expected to render $x^{\star}$
practically asymptotically stable for the closed loop, with the region of
practical convergence shrinking as the horizon $N$ grows. A mathematically
rigorous treatment of this behaviour -- with the precise assumptions and
proofs -- is given in \cite[Chapter~8]{grune2016nonlinear} and the references
therein. Establishing those hypotheses formally for the present
nine-dimensional epidemic model is beyond the scope of this paper; we instead
demonstrate the turnpike behaviour directly and numerically in
Section~\ref{subsec:turnpike}.

We do not impose the turnpike property by assumption; we exhibit it directly. The
open-loop experiments of Section~\ref{subsec:turnpike} display its defining
signature -- trajectories from widely separated initial states collapse onto a
common suppression arc and remain there until a terminal exit -- and the
invariance of that arc to the horizon length confirms that the convergence phase
is horizon-independent. Because the suppression arc is held at the active
capacity boundary, the turnpike here is \emph{constraint-induced}: it is shaped by
the healthcare-capacity barrier, a setting addressed within the stability and
feasibility theory of state-constrained economic MPC without terminal
conditions.\cite{boccia2014stability,faulwasser2018economic} Together with the
inherent feedback of receding-horizon control, this directly demonstrated
turnpike underwrites the practical asymptotic stability observed throughout for
the chosen horizon $N=35$.

\section{Simulation Results and Discussion}\label{sec:results}
We present the numerical validation of the proposed SEIR--V MPC framework. The
simulation evaluates a high-stakes ``Peak Crisis'' scenario, assessing the
controller's ability to manage an aggressive epidemic wave while taking into
account economic continuity.

\subsection{Simulation Setup and Biological Parameters}\label{subsec:setup}
The simulation models a population of $N_{\mathrm{pop}}=60{,}000{,}000$. The
biological parameters governing the compartmental transitions are calibrated to a
highly infectious viral strain\cite{giordano2020modelling} and are listed in full
in Appendix~\ref{app:params}. The incubation rate is $\kappa=0.2$ (mean
incubation $5$~days), with recovery rates $\gamma_a=0.2$ and $\gamma_i=0.1$. The
baseline disease-induced mortality rates are $\delta_i=\delta_p=0.005$, while the
critical-care mortality rate is $\delta_h=0.0124$.

Crucially, the simulation is initialized not from a disease-free state but from a
\textbf{Peak Crisis (1~January 2021)} scenario reflecting a heavily burdened
system:
\begin{itemize}
\item Hospitalized $H_0=22{,}822$.
\item Fatalities $F_0=74{,}159$ cumulative deaths.
\item Initial infections $I_1(0)=11{,}411$ and exposed $E(0)=22{,}822$.
\item The susceptible population $S(0)$ is the remainder of $N_{\mathrm{pop}}$
minus all other compartments, totaling approximately $58.3$~million.
\end{itemize}

\subsection{Aggressive Control Strategy and the ``Deadzone'' Cost
Function}\label{subsec:aggressive}
Unlike conservative elimination strategies, the controller is tuned to execute an
``Aggressive'' or endemic-maintenance policy:\cite{lavine2021immunological} it
accepts a continuous, non-zero level of infection, driving $\beta(t)$ as close to
the fully open state ($\beta_{\max}$) as the health constraint permits and
intervening to prevent hospitalizations from breaching capacity. The control
variables are bounded by physical and logistical realities: the vaccination rate
is capped at $w_{\max}=0.01$ (1\% of the population per day), while the
transmission bounds are $\beta_{\min}=0.17$ (hard lockdown) and $\beta_{\max}=0.50$
(fully open).

The aggressive policy is realized through an asymmetric, priority-driven weight
distribution:
\begin{enumerate}
\item \textbf{Economic Priority ($w_{\mathrm{econ}}=8000$):} the penalty for
restricting the economy is set exceptionally high, forcing the controller to
``fight for every dollar'' of GDP \cite{alvarez2021simple} (a baseline daily
output of approximately \$5.5~billion, consistent with Italy's annual GDP of about
\$2~trillion).\cite{worldbank_gdp_italy}
\item \textbf{Healthcare ``Deadzone'' ($H_{\mathrm{safe}}=25{,}000$):} if
$H(t)\le 25{,}000$ the health cost is identically zero; beyond it a massive
quadratic penalty ($w_{\mathrm{health}}=1000$) activates, normalized by the
remaining margin of $5{,}000$ beds (the gap to the absolute limit
$H_{\max}=30{,}000$), acting as a steep barrier.\cite{rawlings2017model}
\item \textbf{Mortality Tolerance ($w_{\mathrm{death}}=100$):} the daily fatality
weight is deliberately moderate, representing a mathematically necessary
acceptance of endemic mortality to preserve the macro-economy.
\item \textbf{Maximum Vaccination ($w_{\mathrm{vax}}=0.1$):} the vaccination
penalty is negligible, so the optimal policy saturates the vaccination input
($w_1=0.01$).
\end{enumerate}

The weights encode a deliberate \emph{lexicographic} priority -- economy above the
capacity barrier, the barrier above tolerated mortality, mortality above
vaccination effort -- chosen to instantiate the most adversarial ``living with
COVID'' stance against which a public-health controller can be stress-tested.
Their absolute magnitudes are not themselves the object of study: the turnpike
analysis of Section~\ref{subsec:turnpike} shows that the qualitative
suppression-then-reopening (Hammer-and-Dance) structure is invariant across
initial conditions and horizons, so it is a structural property of the
constrained problem rather than an artifact of the particular weights. A
systematic Pareto-front exploration of the health--economy trade-off across
alternative weightings is a natural and worthwhile extension, which we leave to
future work. The monetary figures reported below follow a transparent accounting:
the daily economic loss is modeled as a fixed fraction of baseline daily output
scaled by squared severity, with auxiliary per-dose and per-bed-day costs, as
detailed in Appendix~\ref{app:params}.

\subsection{Analysis of Closed-Loop Trajectories}\label{subsec:closed}
The closed-loop simulation was executed over $149$~days with a prediction horizon
$N=35$ and sampling time $T_s=1$~day.

\begin{figure}[!ht]\centering
\begin{minipage}[t]{0.48\textwidth}\centering
\includegraphics[width=\linewidth]{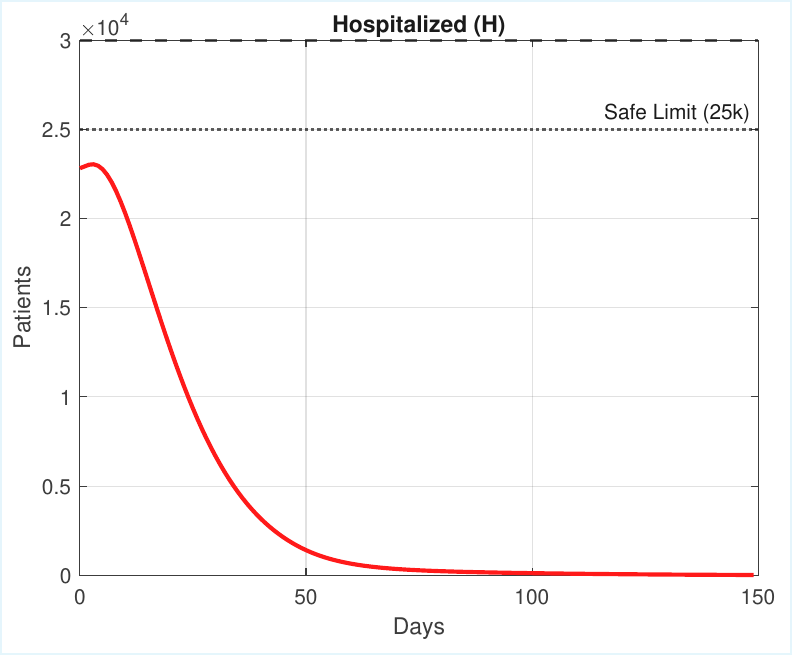}\\[2pt]{\small (a)}
\end{minipage}\hfill
\begin{minipage}[t]{0.48\textwidth}\centering
\includegraphics[width=\linewidth]{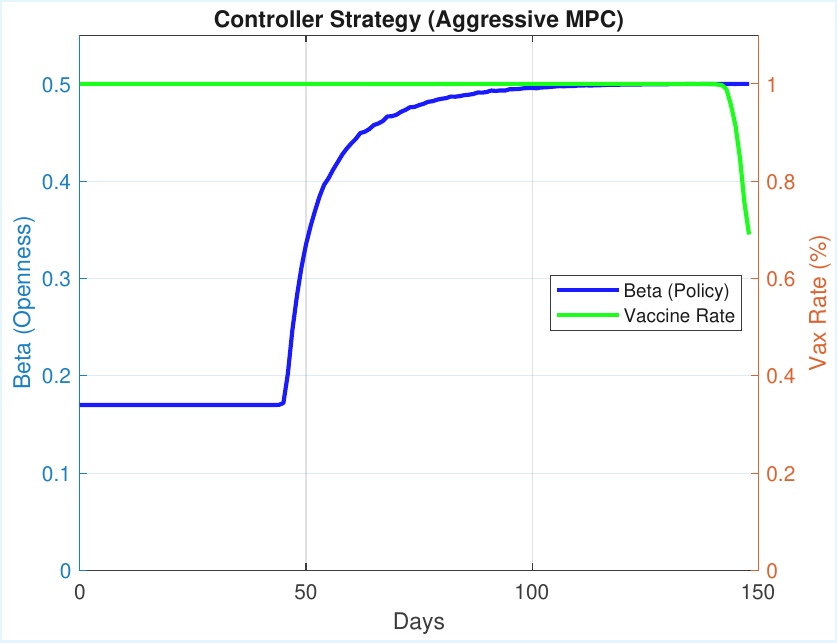}\\[2pt]{\small (b)}
\end{minipage}

\vspace{8pt}
\begin{minipage}[t]{0.48\textwidth}\centering
\includegraphics[width=\linewidth]{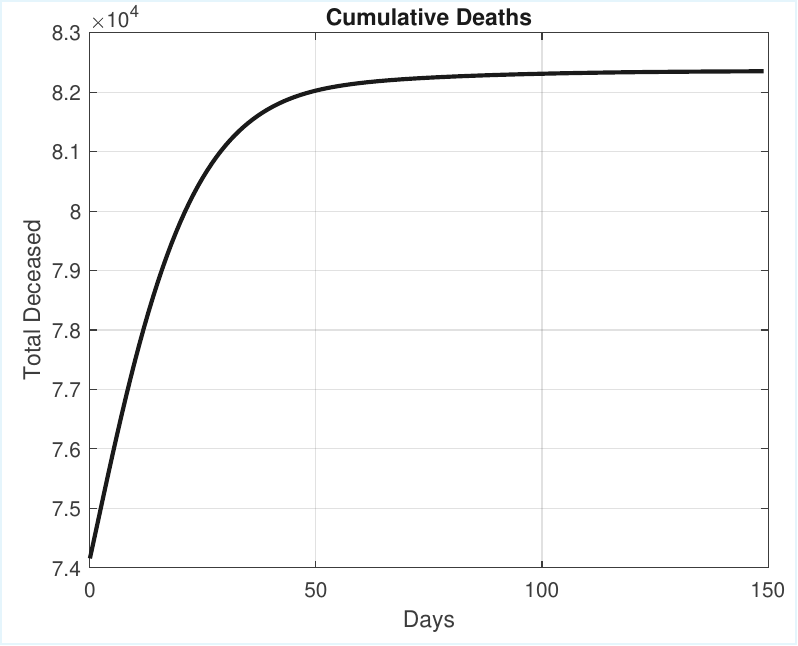}\\[2pt]{\small (c)}
\end{minipage}\hfill
\begin{minipage}[t]{0.48\textwidth}\centering
\includegraphics[width=\linewidth]{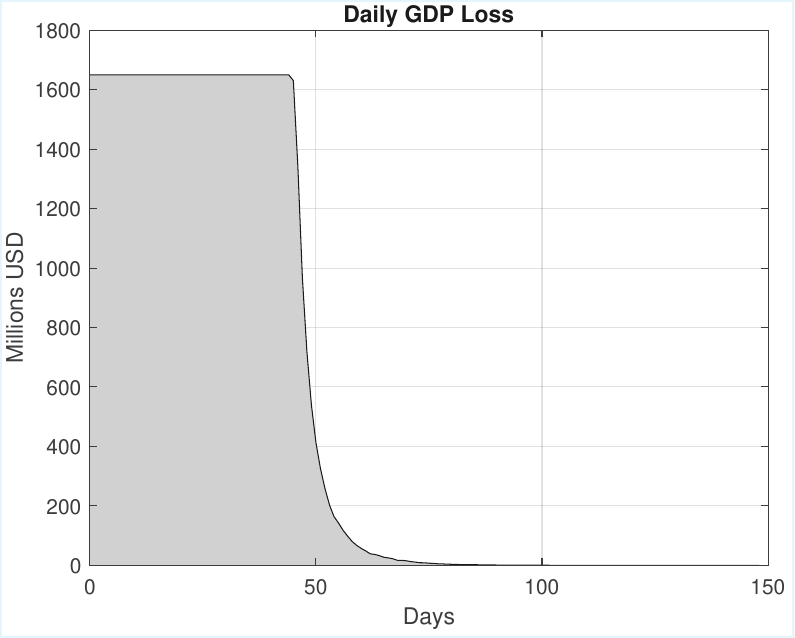}\\[2pt]{\small (d)}
\end{minipage}
\caption{Aggressive Living-with-COVID strategy (closed loop, $N=35$):
(a) hospitalizations $H(t)$ peak at $22{,}859$ and are driven below the
$25{,}000$ deadzone toward zero; (b) the controller holds $\beta=\beta_{\min}$
(the ``Hammer'') for $\sim45$ days, then reopens toward $\beta_{\max}$ (the
``Dance''), with vaccination maximized throughout; (c) cumulative deaths plateau
at $82{,}318$; (d) daily GDP loss, high during suppression and collapsing after
reopening.}
\label{fig:closed}
\end{figure}

Starting from the elevated load $H_0=22{,}822$, dangerously close to the
$25{,}000$ deadzone, the controller immediately applies maximum suppression,
holding the economic openness at its lower bound $\beta=\beta_{\min}=0.17$ (hard
lockdown) for approximately the first $45$~days (Figure~\ref{fig:closed}b). This
is the ``Hammer'': together with vaccination held at its $1\%$ ceiling throughout,
it drives the hospitalized compartment monotonically down from its peak of
$22{,}859$ toward zero (Figure~\ref{fig:closed}a), so that $H$ never breaches the
$25{,}000$ limit. Crucially, the controller does \emph{not} open the economy
prematurely: with the highly transmissible mutant strain ($c_2=1.5$), any early
relaxation would drive $H$ through the capacity barrier and incur the steep
deadzone and fatality penalties, so suppression is optimal despite the dominant
economic weight.

Once the susceptible pool has been sufficiently depleted by vaccination and the
active infection suppressed, the marginal health risk of reopening collapses and
the economic penalty takes over: from around day~$45$ the controller
progressively raises $\beta(t)$ toward $\beta_{\max}=0.50$ -- the ``Dance'' --
completing the suppression-then-reopening arc (Figure~\ref{fig:closed}b).
Cumulative fatalities accrue almost entirely during the early suppression phase,
rising from $74{,}159$ to a plateau of $82{,}318$ by roughly day~$60$
(Figure~\ref{fig:closed}c); correspondingly, the daily economic loss is high
while the lockdown holds and falls to near zero after reopening
(Figure~\ref{fig:closed}d). Notably, this Hammer-and-Dance policy is not
prescribed: it emerges autonomously from the receding-horizon optimization, and
it mirrors exactly the open-loop turnpike behavior analyzed in
Section~\ref{subsec:turnpike}.

\subsection{Sensitivity Analysis of Prediction Horizon}\label{subsec:sensitivity}
\begin{figure}[!ht]\centering
\includegraphics[width=0.72\linewidth]{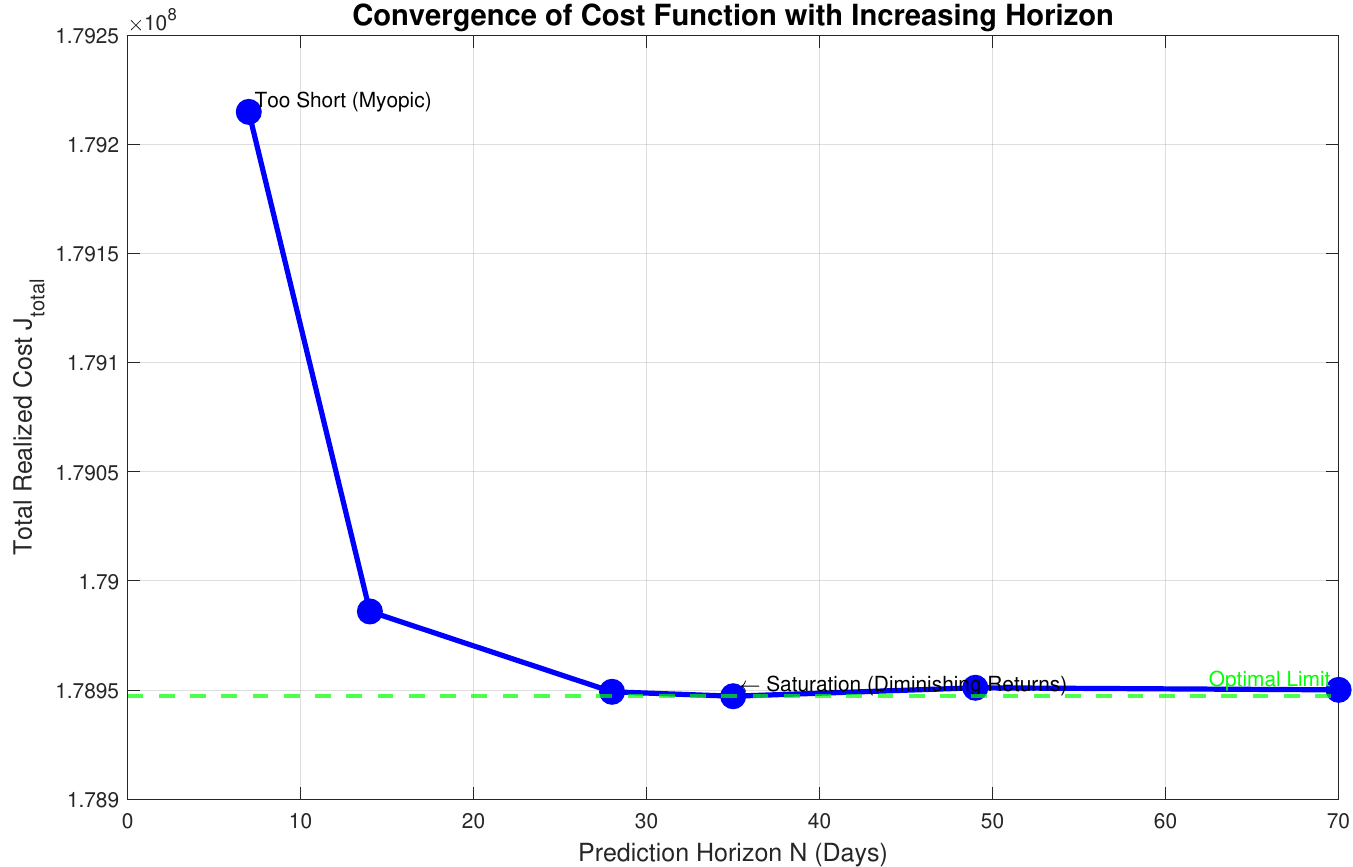}
\caption{Horizon saturation analysis: total realized cost $J_{\mathrm{total}}$
versus prediction horizon $N$. The cost is minimized at $N=35$ and varies by
under $0.15\%$ across the range.}
\label{fig:sat}
\end{figure}

The prediction horizon $N$ dictates the controller's ability to foresee and
preemptively react to constraint boundaries, but increasing $N$ raises NLP
complexity. We evaluated the closed-loop performance over
$N\in\{7,14,28,35,49,70\}$ days, each simulated from the peak-crisis state over
$149$~days, recording the total realized cost $J_{\mathrm{total}}$
(Table~\ref{tab:sat}, Figure~\ref{fig:sat}).

\begin{table}[!ht]\centering
\caption{Realized closed-loop cost versus horizon $N$. The minimum is at $N=35$; the total spread is below $0.15\%$.}
\label{tab:sat}
\begin{tabular}{rrr}
\hline
$N$ (days) & $J_{\mathrm{total}}$ & Change \\
\hline
$7$  & $179{,}214{,}867$ & --- \\
$14$ & $178{,}986{,}027$ & $-0.128\%$ \\
$28$ & $178{,}949{,}225$ & $-0.021\%$ \\
$35$ & $\mathbf{178{,}947{,}193}$ & $-0.001\%$ \\
$49$ & $178{,}951{,}068$ & $+0.002\%$ \\
$70$ & $178{,}950{,}115$ & $-0.001\%$ \\
\hline
\end{tabular}
\end{table}

The results reveal a closed loop that is \emph{robust} to the horizon: the
realized cost varies by under $0.15\%$ across the entire range, with a
well-defined minimum at $N=35$. Three regimes are nonetheless discernible.
\begin{enumerate}
\item \textbf{The Myopic Regime ($N\le 14$).} The shortest horizons incur the
highest cost. The penalty is modest, however, because the elevated initial load
makes the imminent hospitalization pressure evident even within a one-week
window: every horizon therefore initiates the same suppression, and the
differences arise only in the fine timing of the subsequent reopening (exit)
phase.
\item \textbf{The Transition Regime ($14<N<35$).} As the horizon expands, the
reopening is timed more efficiently and the realized cost
falls.\cite{grune2013economic}
\item \textbf{The Saturation Zone ($N\ge 35$).} Beyond $N=35$ the curve
saturates, attaining its minimum at $N=35$ with sub-$0.1\%$ non-monotone
fluctuations thereafter. This non-monotonicity is characteristic of Economic MPC
on non-convex systems without terminal constraints, where the growing NLP
occasionally admits slightly suboptimal local minima.
\end{enumerate}
Based on this saturation test, $N=35$ days is the cost-optimal horizon --
sufficient to time the Hammer-and-Dance exit correctly while keeping the NLP
tractable, avoiding both the sub-optimality of myopic control and the latency and
increased numerical sensitivity of excessively long windows.

\subsection{Turnpike and Stability Analysis}\label{subsec:turnpike}
We now exhibit, in open loop, the turnpike property discussed in
Section~\ref{subsec:feas} (Remark~1).
The turnpike property implies that for a sufficiently long horizon $N$ the
optimal trajectories rapidly converge to and spend most of their time on a
specific path (the ``turnpike''), deviating only near the horizon end due to the
finite-time nature of the optimization -- the ``exit phenomenon.'' We verify this
through two open-loop experiments.

\subsubsection{Experiment 1: Sensitivity to Initial Conditions}
Fixing the horizon and evaluating the open-loop response across three infection
scenarios -- Low ($0.1\times$), Baseline ($1.0\times$) and High ($5.0\times$) --
yields two observations.

\begin{figure}[!ht]\centering
\begin{minipage}[t]{0.48\textwidth}\centering
\includegraphics[width=\linewidth]{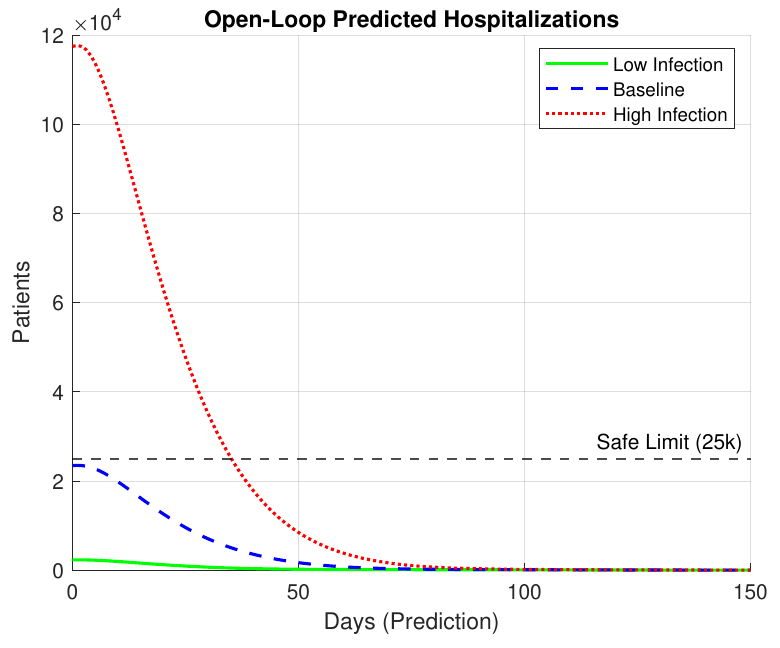}\\[2pt]{\small (a)}
\end{minipage}\hfill
\begin{minipage}[t]{0.48\textwidth}\centering
\includegraphics[width=\linewidth]{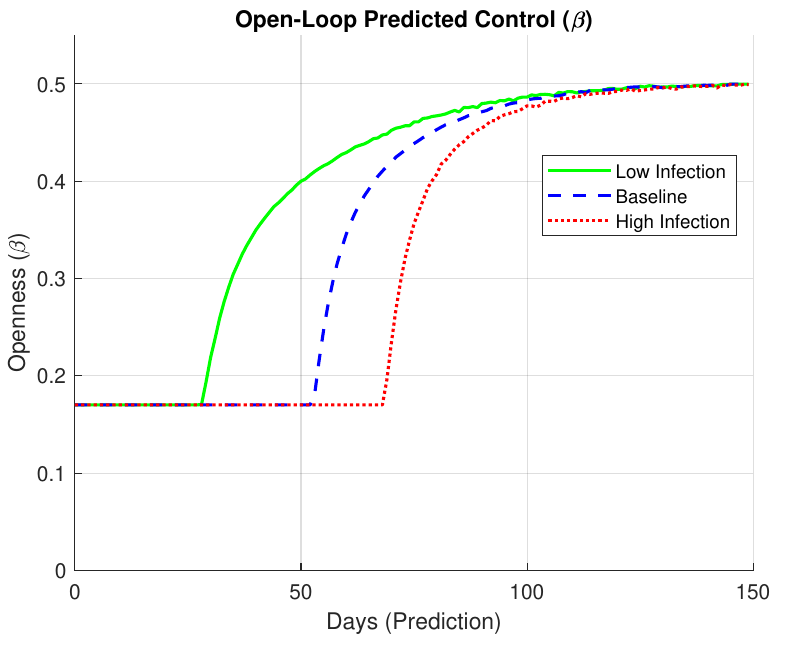}\\[2pt]{\small (b)}
\end{minipage}
\caption{Turnpike property under varying initial conditions: (a) hospitalizations
from $0.1\times$, $1.0\times$ and $5.0\times$ baseline infection converge to a
common optimal boundary path; (b) the openness $\beta(t)$ anchors at the
suppression bound before the progressive exit, regardless of the initial state.}
\label{fig:turnpikeIC}
\end{figure}

\textbf{State Convergence (Hospitalizations).} All three hospitalization
trajectories ultimately converge toward zero (the disease-free equilibrium). In
the High scenario, the immense penalty of violating the soft constraint forces a
prolonged suppression that drives hospitalizations rapidly down through the
$25{,}000$ safe limit (Figure~\ref{fig:turnpikeIC}a). In all scenarios the
optimizer determines that driving active infections near zero is the most
cost-effective route to long-term economic freedom.

\textbf{Control Convergence (Openness $\beta$).} Regardless of the initial state,
the optimal policy anchors to the turnpike at maximum lockdown ($\beta=0.17$) to
crush the infection curve, then exhibits the characteristic exit phenomenon,
progressively stepping up openness ($\beta\to 0.50$) as the threat clears
(Figure~\ref{fig:turnpikeIC}b). The time spent on the suppression turnpike scales
with the initial severity: reading from Figure~\ref{fig:turnpikeIC}b, it is
approximately $60$, $45$ and $20$ days for the High, Baseline and Low scenarios,
respectively.

\subsubsection{Experiment 2: Prediction Horizon and the Exit Phenomenon}
\begin{figure}[!ht]\centering
\begin{minipage}[t]{0.48\textwidth}\centering
\includegraphics[width=\linewidth]{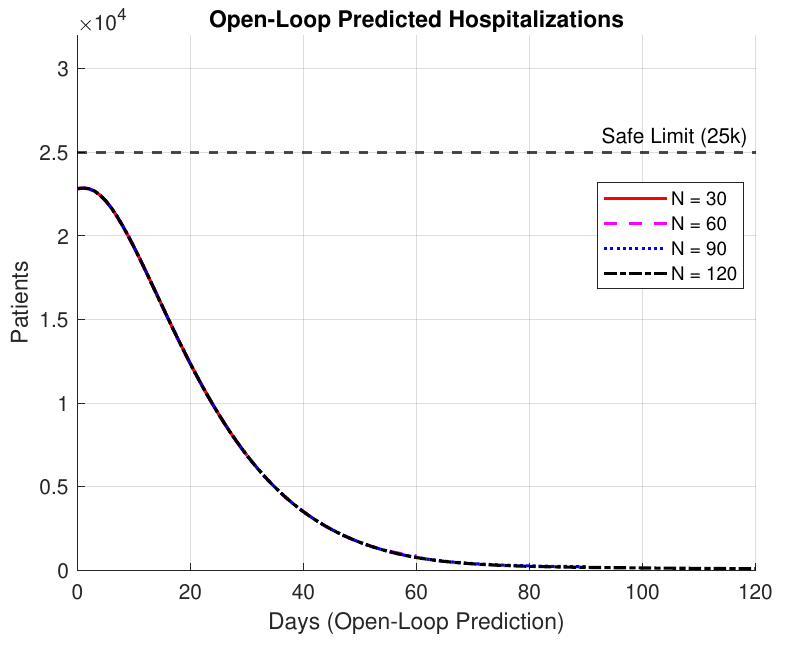}\\[2pt]{\small (a)}
\end{minipage}\hfill
\begin{minipage}[t]{0.48\textwidth}\centering
\includegraphics[width=\linewidth]{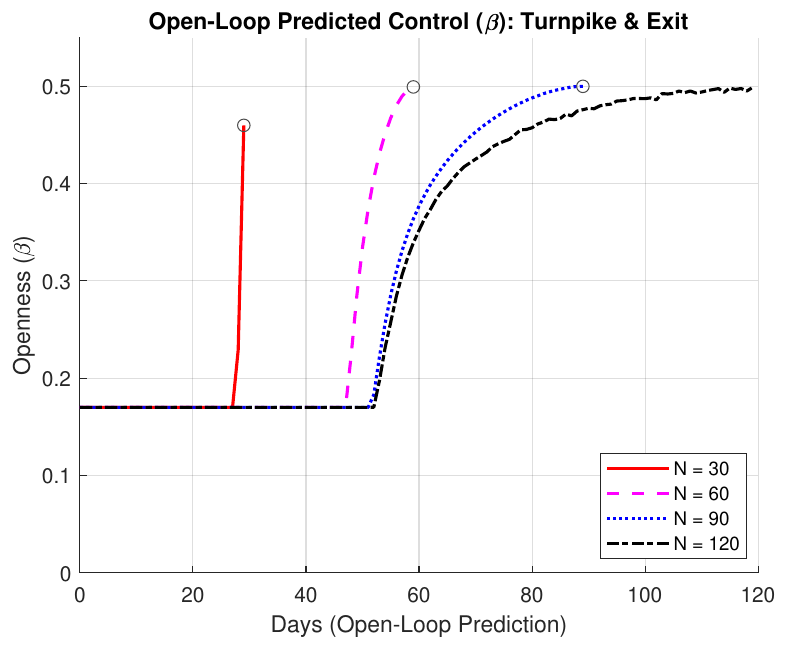}\\[2pt]{\small (b)}
\end{minipage}
\caption{Open-loop turnpike analysis (varying $N$): (a) states -- longer horizons
spend a greater proportion of time on the turnpike, shorter horizons diverge
earlier; (b) control -- the characteristic exit phenomenon as the optimizer
becomes greedy near the horizon end.}
\label{fig:turnpikeN}
\end{figure}

Holding the baseline initial conditions fixed and varying $N\in\{30,60,90,120\}$
isolates the stationary turnpike phase from the transient exit phase. In the
state plot all trajectories overlap for an identical initial segment
(Figure~\ref{fig:turnpikeN}a); in the control plot the trajectories for longer
horizons remain on the turnpike at $\beta\approx 0.17$ before each, as it
approaches its horizon, exhibits a ``greedy exit'' in which the economy is opened
prematurely because the finite-horizon optimizer no longer accounts for cost
beyond day $N$ (Figure~\ref{fig:turnpikeN}b).

\subsubsection{Implications for Stability}
Before drawing conclusions we note a minor residual variation in the openness
trajectories for the longest horizons near the horizon end
(Figure~\ref{fig:turnpikeN}b), together with a very slight rise of the realized
cost beyond its minimum for the largest horizons in Figure~\ref{fig:sat}. These
are numerical artifacts of the SQP solver applied to the growing, non-convex
NLP -- a consequence of its finite convergence tolerance and its occasional
settling on marginally different local minima as the problem dimension
increases -- rather than genuine features of the optimal policy. They do not
affect the qualitative Hammer-and-Dance structure, the location of the cost
minimum at $N=35$, or the turnpike behaviour discussed below.
These open-loop analyses provide the turnpike signature predicted by
Definition~1; as discussed in Section~\ref{subsec:feas}, it is this
constraint-induced turnpike that underlies the practical asymptotic stability of
the optimal operating point.\cite{boccia2014stability} The invariance of the
initial trajectory segments to the horizon length $N$ ensures that the
receding-horizon (closed-loop) implementation remains stable and optimal,
provided $N$ exceeds the initial convergence phase (verified as $N\ge 35$ in the
sensitivity analysis). The numerical evidence thus indicates that the capacity
barrier steers the population toward the optimal economic--health operating point
without loss of recursive feasibility, even under extreme initial-state
perturbations.

\section{Conclusion}\label{sec:conclusion}
The COVID-19 pandemic underscored the fragility of healthcare systems and the
socio-economic cost of broad-spectrum NPIs. We addressed the challenge of
designing an optimal control strategy that balances saving lives with preserving
economic continuity, formulating pandemic management within a nonlinear
\textbf{Model Predictive Control} framework applied to a high-fidelity SEIR--V
model. Our investigation yields four conclusions.
\begin{enumerate}
\item \textbf{Efficacy of Dynamic Feedback.} The Receding Horizon Control scheme
manages the epidemic trajectory effectively. By continuously adjusting $\beta(t)$
in response to real-time data, the closed-loop controller navigates the competing
objectives and holds hospitalizations below the deadzone throughout, without
violating the healthcare capacity constraint.\cite{anderson2020will}
\item \textbf{Optimality of the ``Hammer and Dance''.} The celebrated ``Hammer
and Dance'' strategy \cite{pueyo2020coronavirus,morris2021optimal} emerges naturally
as the mathematical optimum -- in both the closed loop and the open-loop turnpike
analysis. The solver, not the modeler, identifies that an aggressive suppression
phase (``the Hammer'') is required to reset the dynamics, followed by a modulated
reopening (``the Dance''), confirming that early, decisive action is economically
optimal in the long run.
\item \textbf{The Horizon Effect and the Exit Phenomenon.} The closed loop is
robust to the prediction horizon, with the realized cost varying by under
$0.15\%$ and minimized at $N=35$ days; the residual sub-optimality of short
horizons is attributable to the \textbf{exit phenomenon}, in which the optimizer
becomes greedy near the end of a short window. Effective governance therefore
requires a planning window exceeding the combined incubation and generation
periods of the pathogen.
\item \textbf{Structural Stability via the Turnpike Property.} We characterized
the optimal operating point as the disease-free equilibrium -- shown numerically
to be an open-loop \emph{unstable} steady state when $R_0(\beta_{\max})>1$ -- and
established its practical asymptotic stability under the economic MPC feedback
through the turnpike property (Definition~1), demonstrated directly by the
open-loop experiments. This turnpike is \emph{constraint-induced} -- enforced by
the capacity barrier rather than by a strictly convex economic cost: regardless
of the initial infection level, the optimal trajectories anchor to a unique
suppression turnpike before progressively reopening, without loss of recursive
feasibility.\cite{boccia2014stability}
\end{enumerate}
In summary, the conflict between public health and economic stability need not be
treated as a fixed trade-off: it can instead be posed and solved as a dynamic
optimization problem. The present analysis is
deterministic; the robustness established here is that of receding-horizon
feedback to initial-state perturbations and to the residual mismatch absorbed by
re-optimization. A formal treatment of parametric and measurement uncertainty --
via stochastic or tube-based MPC, validated by Monte~Carlo assessment of the
controller against a perturbed plant -- is the subject of ongoing work. By
leveraging the predictive power of MPC and the stability of the Turnpike
Property, governments can navigate this trade-off with minimized loss of life and
maximized economic resilience.

\appendix
\section{Model Parameters and Implementation Details}\label{app:params}
For reproducibility, this appendix collects the calibrated parameters, the
initial state, the control configuration, and the economic accounting used to
generate every result in Section~\ref{sec:results}. The biological parameters are
those identified in the companion calibration study;\cite{melhani2026nine} the
state-estimation layer (Extended Kalman Filter design and tuning) is documented in
the companion estimation study,\cite{melhani2026extended} while the closed-loop
results here use full-state feedback.

\begin{table}[!ht]\centering
\caption{Calibrated model parameters (corrected model, $c_P=c_2=1.5$).}
\label{tab:appparams}
\footnotesize
\begin{tabular}{lll}
\hline
Symbol & Value & Description \\
\hline
$\Lambda$    & $100$~day$^{-1}$    & Recruitment (births) \\
$\mu$        & $3.535\times10^{-5}$ & Natural mortality rate \\
$\kappa$     & $0.20$              & Incubation rate \\
$\gamma_a$   & $0.20$              & Hospitalization rate \\
$\gamma_i$   & $0.10$              & Community recovery rate \\
$\gamma_r$   & $0.10$              & In-hospital recovery rate \\
$\rho_1$     & $0.58$              & Fraction to original strain \\
$\rho_2$     & $0.001$             & Fraction to super-spreaders \\
$m$          & $0.005$             & Mutation rate $I_1\!\to\!I_2$ \\
$r_1,r_2$    & $0.05$              & Direct recovery rates \\
$\delta_i$   & $0.005$             & Mortality, infectious \\
$\delta_p$   & $0.005$             & Mortality, super-spreaders \\
$\delta_h$   & $0.0124$            & Mortality, hospitalized \\
$\delta$     & $0.001$             & Waning natural immunity \\
$\psi$       & $0.002$             & Waning vaccine immunity \\
$\sigma$     & $0.80$              & Vaccine protection factor \\
$c_P,c_2$    & $1.5$               & Transmissibility multipliers \\
$N_{\mathrm{pop}}$ & $6\times10^{7}$ & Population size \\
\hline
\end{tabular}
\end{table}

\noindent\textbf{Initial state (1~January 2021).}
$x_0=[\,S_0,E_0,V_0,I_{1,0},I_{2,0},P_0,H_0,R_0,F_0\,]^{\top}$ with
$S_0\approx5.833\times10^{7}$, $E_0=22{,}822$, $V_0=32{,}812$,
$I_{1,0}=11{,}411$, $I_{2,0}=1{,}141$, $P_0=571$, $H_0=22{,}822$,
$R_0=1{,}500{,}000$, $F_0=74{,}159$, with $S_0$ taken as the residual of
$N_{\mathrm{pop}}$.

\noindent\textbf{Control configuration.} Bounds $\beta\in[0.17,0.50]$,
$w_1\in[0,0.01]$; capacity $H_{\mathrm{safe}}=25{,}000$, $H_{\max}=30{,}000$,
$H_{\mathrm{scale}}=5{,}000$; cost weights $w_{\mathrm{econ}}=8000$,
$w_{\mathrm{health}}=1000$, $w_{\mathrm{death}}=100$, $w_{\mathrm{vax}}=0.1$,
$w_{\mathrm{smooth}}=1000$; horizon $N=35$, sampling $T_s=1$~day, simulation
length $149$~days.

\noindent\textbf{Economic accounting.} With baseline daily output
$G=\$5.5\times10^{9}$ and severity
$s_k=(\beta_{\max}-\beta_k)/(\beta_{\max}-\beta_{\min})$, the modeled daily
economic loss is $\mathcal{L}(k)=0.30\,G\,s_k^{2}$. Reported auxiliary monetary
totals use \$20 per vaccine dose and \$2{,}000 per hospital bed-day. These
conversions affect only the reported dollar figures, not the control law, which
depends solely on the dimensionless stage cost $\ell$.

\section*{Acknowledgments}
The doctoral scholarship of L.~R.~Melhani is funded by InEmbryo S.r.l.s.\ through
the University of Palermo. The authors thank the Italian Civil Protection
Department for the publicly available epidemiological data.

\section*{Conflict of Interest}
The authors declare no potential conflict of interest.

\section*{Data Availability Statement}
The Italian Civil Protection COVID-19 data used for calibration are publicly
available.


\end{document}